\documentclass[11pt]{amsart}
\usepackage[margin=1in]{geometry}
\usepackage[T1]{fontenc}
\usepackage{lmodern}
\usepackage{amsmath,amssymb,mathtools}
\usepackage{microtype}
\usepackage{tikz,placeins}
\usetikzlibrary{arrows.meta}
\usepackage[hidelinks]{hyperref}

\DeclareMathOperator{\conv}{conv}
\DeclareMathOperator{\SC}{SC}
\DeclareMathOperator{\rec}{rec}
\DeclareMathOperator{\dist}{dist}
\DeclareMathOperator{\linspan}{span}
\DeclareMathOperator{\cone}{cone}
\newcommand{\R}{\mathbb R}
\newcommand{\Z}{\mathbb Z}
\definecolor{inkblue}{HTML}{253B53}
\definecolor{cutteal}{HTML}{007F82}

\theoremstyle{plain}
\newtheorem{theorem}{Theorem}
\newtheorem{lemma}{Lemma}
\newtheorem{proposition}{Proposition}

\title{A Polynomial Bound on Split Rank}
\author{Avinash Bhardwaj}
\address{Department of Industrial Engineering and Operations Research, Indian Institute of Technology Bombay, Mumbai 400076, India}
\email{abhardwaj@iitb.ac.in}
\author{Vishnu Narayanan}
\address{Department of Industrial Engineering and Operations Research, Indian Institute of Technology Bombay, Mumbai 400076, India}
\email{vishnu@iitb.ac.in}
\date{}
\hypersetup{pdftitle={A Polynomial Bound on Split Rank},pdfauthor={}}

\begin{document}
\begin{abstract}
We prove that the split rank of a rational pure-integer polyhedron is polynomially bounded in its input length, resolving a conjecture of Basu, Conforti, Di Summa, and Jiang. If the defining system has dense binary encoding length $L$, then $O((L+1)^9)$ split rounds suffice to reach the integer hull in arbitrary dimension.
\end{abstract}
\maketitle

\section{Introduction}
Let $P=\{y\in\R^n:a_i^\top y\le b_i,\ 1\le i\le m\}$ be a rational polyhedron with integer hull $P_I=\conv(P\cap\Z^n)$. Throughout, we assume without loss of generality that the defining inequalities of $P$ and $P_I$ and the target inequality $g\le0$ have integer coefficients. A \emph{split} $(\pi,k)$, where $\pi\in\Z^n\setminus\{0\}$ and $k\in\Z$, has two branches, $\pi^\top x\le k$ and $\pi^\top x\ge k+1$. Every integer point satisfies one of them. An inequality valid on both branches in $P$ is a \emph{split cut} for $P$. One \emph{split round} imposes all split cuts for the current relaxation. Starting from $\SC^0(P)=P$, write $\SC^j(P)$ for the relaxation after $j$ rounds. The \emph{split rank} of $P$ is the least integer $j\ge0$ such that $\SC^j(P)=P_I$. Schrijver proved that finitely many Chv\'atal--Gomory rounds suffice to reach $P_I$~\cite[Theorem~2]{Schrijver80}. Since every Chv\'atal--Gomory cut is a split cut, the split rank is finite. When $P_I=\varnothing$, Cook, Coullard, and Tur\'an~\cite[Theorem~$1'$]{CCT87} proved that a number of Chv\'atal cuts depending only on the dimension suffices to reduce $P$ to the empty set. For mixed-integer polyhedra, even one continuous variable can give infinite split rank~\cite[Example~2]{CKS}.

In the pure-integer case, every rational polyhedron in the plane has split rank at most two~\cite[Remark~3.15]{BCDJplane}. In dimension three, the rank can be arbitrarily large even when the integer hull is fixed. To see this, let $T$ be the triangle with vertices $(0,0,0)$, $(2,0,0)$, and $(0,2,0)$. For a positive integer $H$, let $T_H$ be the tetrahedron obtained by adjoining the apex $(1/2,1/2,H)$. Every point of $T_H$ at positive height satisfies $x_1>0$, $x_2>0$, and $x_1+x_2<2$. No integer pair satisfies these inequalities. Since the vertices of $T$ are integral, $T$ is the integer hull of $T_H$. Conforti et al. showed that deriving $x_3\le0$ on $T_H$ requires $\Omega(\log H)$ split rounds~\cite[Lemma~19]{CDDFG}. The four inequalities defining $T_H$ take $O(\log(H+1))$ bits in total, so this family gives a lower bound linear in the input length. More generally, Conforti et al. characterize the integer hulls that admit rational relaxations of arbitrarily large split rank~\cite[Theorem~1]{CDDF}.

Basu, Conforti, Di Summa, and Jiang~\cite[Section~4]{BCDJComplexity} conjectured that the split rank of a rational polyhedron is bounded by a polynomial in its input length, even when the dimension is part of the input. In~\cite{BN26}, the authors established the conjecture in fixed dimension by constructing a split proof of polynomial total encoding length. The polynomial bound in~\cite{BN26} depends on the dimension. The following theorem bounds the number of split rank by a single polynomial when the dimension is part of the input.

\begin{theorem}\label{thm:generalrank}
Let $P\subseteq\R^n$ be a rational pure-integer polyhedron given by a system of linear inequalities of dense binary encoding length $L$. The split rank of $P$ is $O((L+1)^9)$.
\end{theorem}

A split round imposes all split cuts for the current relaxation. A split proof adds one cut at a time, and its length is the number of cuts it uses. For a polytope in $[0,1]^n$, taking the convex hull over the two branches of each coordinate split $(e_i,0)$ in succession gives the integer hull in at most $n$ rounds~\cite{Balas98}. This bound is tight~\cite{CL02}. Even a polytope in the cube defined by $O(n)$ inequalities can have empty integer hull and split rank $\Omega(n)$~\cite{DD24}. Yet some $0$--$1$ systems require exponentially many individual split cuts to prove infeasibility~\cite{Dash10}. In Section~\ref{sec:encoding} we show that the round bound in Theorem~\ref{thm:generalrank} still holds when each cut and its split have encoding length polynomial in $L$.

\subsection*{Idea of the proof}
For the tetrahedron $T_H$, deriving $x_3\le0$ suffices to reach the integer hull. More generally, suppose $P_I\ne\varnothing$, and let $g(x)=c^\top x+c_0\le0$ be a defining inequality of $P_I$ that is violated somewhere on $P$. We call $g(x)$ the height of $x$. Since $P_I\ne\varnothing$, $P$ and $P_I$ have the same recession cone~\cite{Schrijver86}. Validity of $g\le0$ on $P_I$ therefore gives $c^\top r\le0$ for every $r\in\rec(P)$. Hence $g$ is bounded above on $P$.

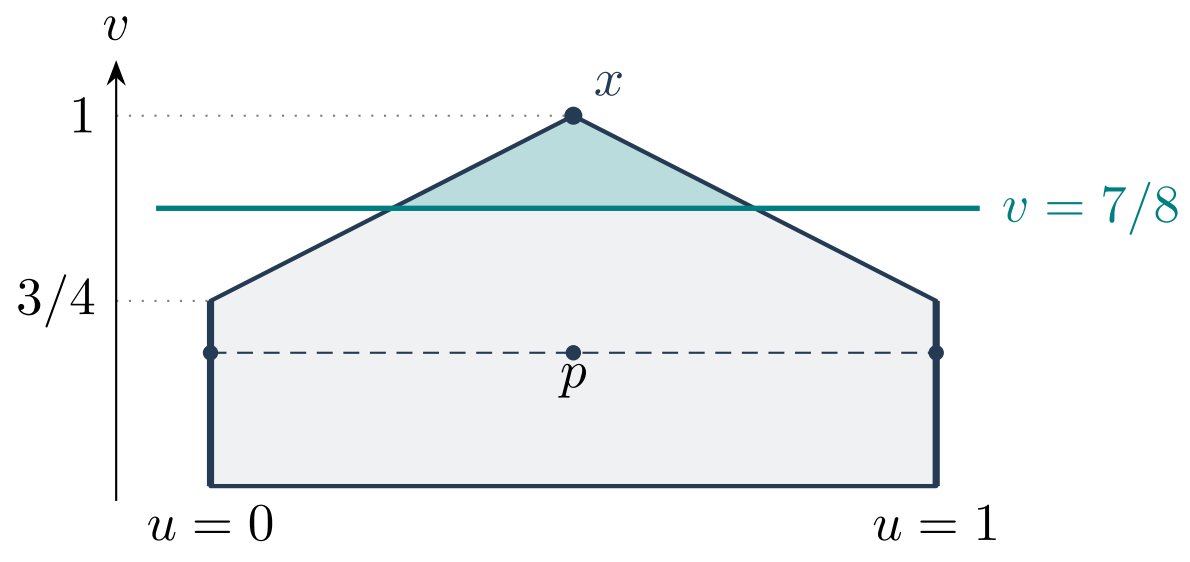
\begin{figure}[!htb]
\centering
\begin{tikzpicture}[x=4.8cm,y=4.9cm,>=Stealth,font=\small,line join=round]
  \fill[inkblue!7] (0,.50)--(0,.75)--(.5,1)--(1,.75)--(1,.50)--cycle;
  \fill[cutteal!27] (.25,.875)--(.5,1)--(.75,.875)--cycle;
  \draw[inkblue,thick] (0,.50)--(0,.75)--(.5,1)--(1,.75)--(1,.50)--cycle;
  \draw[inkblue,line width=1.35pt] (0,.50)--(0,.75) (1,.50)--(1,.75);
  \draw[inkblue,densely dashed] (0,.68)--(1,.68);
  \foreach \u in {0,.5,1}{\fill[inkblue] (\u,.68) circle(1.45pt);}
  \node[below,inner sep=2pt] at (.5,.68) {$p$};
  \draw[cutteal,line width=1pt] (-.075,.875)--(1.06,.875) node[right] {$v=7/8$};
  \draw[gray,dotted] (-.13,1)--(.49,1) (-.13,.75)--(0,.75);
  \draw[->] (-.13,.48)--(-.13,1.075) node[above] {$v$};
  \node[left] at (-.13,1) {$1$};
  \node[left] at (-.13,.75) {$3/4$};
  \node[below] at (0,.50) {$u=0$};
  \node[below] at (1,.50) {$u=1$};
  \fill[inkblue] (.5,1) circle(1.7pt) node[above right] {$x$};
\end{tikzpicture}
\caption{The convex function $F(u,v)=v+(u-\tfrac12)^2$ is at most one at the vertices, hence throughout the polygon. Its tangent at $x=(\tfrac12,1)$ is $T(u,v)=v$. On both branches of $(e_1,0)$, $F-T=\tfrac14$, so the raised tangent $T+\tfrac18\le1$ cuts off the cap $v>\tfrac78$ containing $x$. The midpoint $p$ survives every cut from this split because it lies between the two marked branch points.}\label{fig:tangent-cut}
\end{figure}
\FloatBarrier

Let $H=\lceil\max_{x\in P}g(x)\rceil$. We show that polynomially many split rounds suffice to halve the upper bound on $g$. After $1+\lfloor\log_2 H\rfloor$ halvings, $g(x)<1$ for every $x$ in the current relaxation, so the branch $g\ge1$ contains no point of that relaxation. Thus the split $(c,-c_0)$ gives $g\le0$ in one further round. Section~\ref{sec:rankbound} bounds both the number of halvings and the number of rounds per halving uniformly over the defining inequalities of $P_I$. Thus all these inequalities hold after the same polynomial number of rounds.

To see how one halving works, consider a bounded $P$. Fix a point $x$ of height greater than $H/4$. Around $x$, we construct an ellipsoid small enough that every point in it still has positive height and satisfies $a_i^\top y<b_i+1$ for every defining inequality $a_i^\top y\le b_i$ of $P$. An integer point in the ellipsoid would therefore satisfy every defining inequality of $P$ and violate $g\le0$. Thus the ellipsoid contains no integer point. Banaszczyk's theorem~\cite{Ban} then gives a split whose two branches stay a positive distance from $x$, measured relative to the ellipsoid. The lower bound on this distance depends only on the dimension.

\begin{figure}[!htb]
\centering
\begin{tikzpicture}[
  box/.style={draw=inkblue,rounded corners=1.6pt,align=center,inner xsep=6pt,inner ysep=3pt,font=\small},
  arr/.style={-{Stealth[length=3.6pt]},draw=inkblue,semithick},
  ann/.style={font=\footnotesize,align=center,fill=white,inner sep=1.5pt}]
  \node[box,fill=inkblue!5] (distance) at (0,0)
     {At $g(x)>H/4$, use the lattice-free ball in $D_x\R^n$\\to find a split satisfying~\eqref{rank:eq:sep}.};
  \node[box] (cut) at (0,-1.16)
     {Use the tangent gap to cut off each point violating\\the next bound on the convex function.};
  \node[box,fill=cutteal!9] (half) at (0,-2.42)
     {Repeat split rounds until $g\le H/2$.};
  \node[box] (row) at (0,-3.53)
     {Apply the split on $g$ to obtain $g\le0$.};
  \node[box,fill=inkblue!5] (hull) at (0,-4.64)
     {All defining inequalities hold on the same closure.};
  \draw[arr] (distance)--(cut);
  \draw[arr] (cut)--(half);
  \draw[arr] (half.east)--(4.60,-2.42)--(4.60,0)--(distance.east);
  \node[ann,anchor=west] at (4.65,-1.22) {$H/2\ge1$\\set $H\gets H/2$};
  \draw[arr] (half)--node[ann,right] {$H/2<1$} (row);
  \draw[arr] (row)--node[ann,right] {take the maximum round count over the rows} (hull);
\end{tikzpicture}
\caption{Proof of Theorem~\ref{thm:generalrank} for $P_I\ne\varnothing$. Section~\ref{sec:rounds} proves the halving step. Section~\ref{sec:rankbound} bounds both the number of halvings and the rounds per halving uniformly over the defining rows. For nonempty $P$ with $P_I=\varnothing$, one halving with $g\equiv1$ and $H=1$ suffices.}\label{fig:rank-roadmap}
\end{figure}
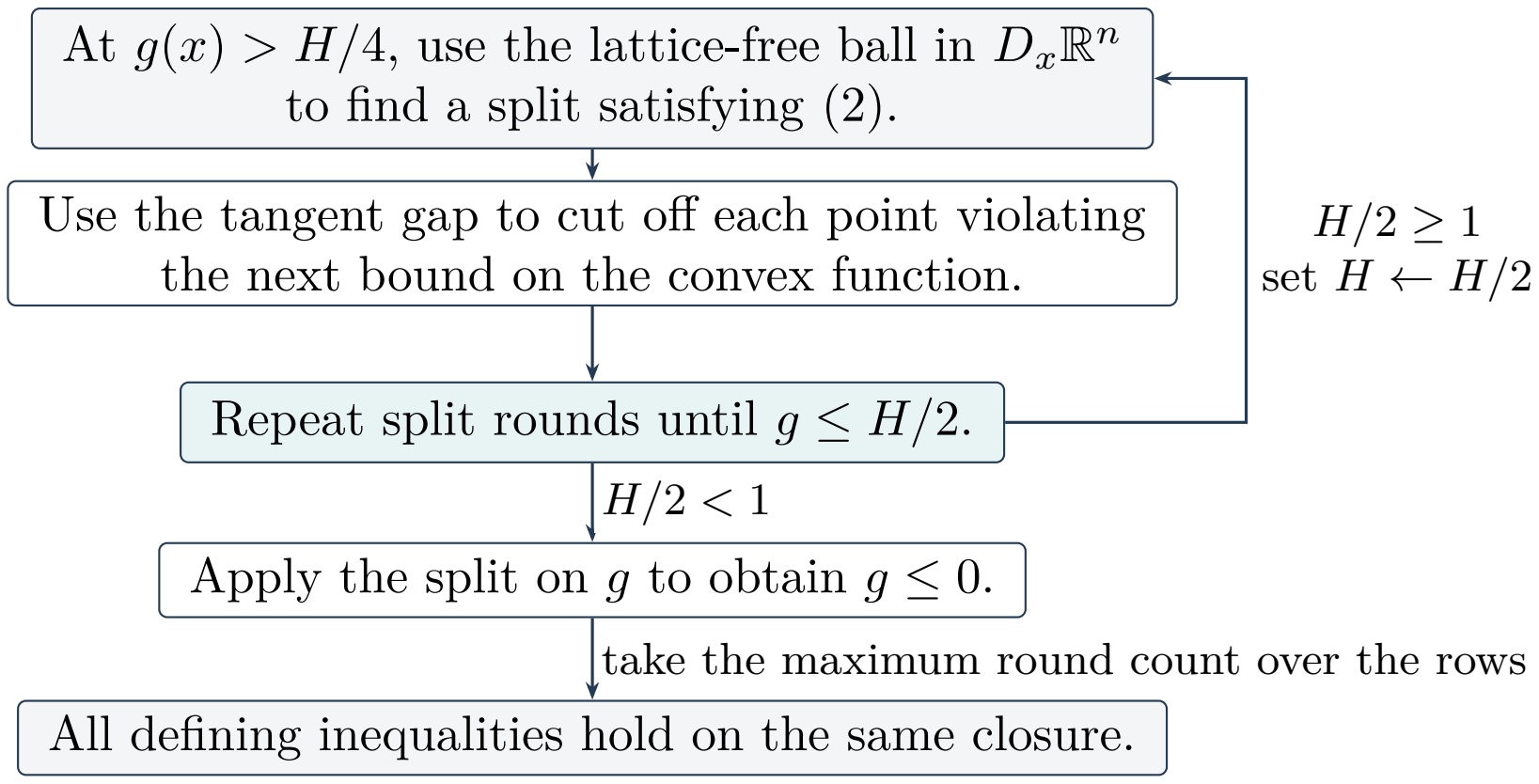

A point between the branches can still belong to their convex hull. To obtain a cut, we use a convex function of the height and the slacks of the defining inequalities. Its tangent at $x$ agrees with the function at $x$. On either branch at nonnegative height, the distance bound gives a gap between the tangent and the function. This leaves room to raise the tangent and obtain an affine inequality that excludes $x$ whenever it violates the bound for the next round. Figure~\ref{fig:tangent-cut} illustrates this step. At heights at least $H/2$, the upper bound on the function decreases in each round until it falls below a fixed lower bound. By then every remaining point has height less than $H/2$. Figure~\ref{fig:rank-roadmap} shows how repeating this halving gives the rank bound.

\FloatBarrier

\section{Halving the maximum violation}\label{sec:rounds}

Let $g(y)=c^\top y+c_0$ be bounded above on $P$, with $g\le0$ on $P_I$. Suppose the current relaxation $Q$ satisfies $P_I\subseteq Q\subseteq P$ and $g\le H$ on $Q$, where $H>0$. We seek to derive $g\le H/2$ by split cuts. If this bound already holds on $Q$, no further rounds are needed. Otherwise, $K=P\cap\{g\ge0\}$ is nonempty. For each defining inequality $a_i^\top y\le b_i$ of $P$, let $s_i(y)=b_i+1-a_i^\top y$ be its shifted slack. Then $s_i\ge1$ on $P$, so $\log s_i$ is defined even where the inequality is tight. For $z\in\Z^n$, $s_i(z)$ is an integer, so $s_i(z)>0$ is equivalent to $a_i^\top z\le b_i$.

\subsection{Distance to the split branches}\label{sec:distance}
To bound the variation of $\log s_i$ on $K$, we also need an upper bound on $s_i$. When $P$ is bounded, every shifted slack has a finite maximum on $K$. For an unbounded example, take $P=[0,\tfrac12]\times[0,\infty)$ and $g(y)=y_1$. Here $K=P$, and the shifted slack $1+y_2$ of the row $y_2\ge0$ is unbounded. We can omit this row because $0\le y_1\le\tfrac12$ already forces $y_1=0$ at every integer point satisfying the remaining inequalities. Lemma~\ref{lem:rows} shows that we can always retain a set of rows whose shifted slacks are bounded on $K$ and whose inequalities exclude every integer point with $g>0$.

\begin{lemma}\label{lem:rows}
There is a set $J\subseteq\{1,\ldots,m\}$ such that $B_i=\max_{y\in K}s_i(y)$ is finite for each $i\in J$, and every $z\in\Z^n$ satisfying $a_i^\top z\le b_i$ for all $i\in J$ has $g(z)\le0$.
\end{lemma}
\begin{proof}
Let $R=\rec(K)$, the set of directions $r$ such that $y+tr\in K$ for every $y\in K$ and $t\ge0$. Along each such direction, $a_i^\top r\le0$ for every $i$, since $a_i^\top(y+tr)\le b_i$ for all $t\ge0$. Moreover, $g(y+tr)=g(y)+t c^\top r$ stays between $0$ and $\max_P g$, so $c^\top r=0$. Since $R$ is a rational polyhedral cone, it has a finite set $E$ of integer generators. Set $r_0=\sum_{e\in E}e$, taking $r_0=0$ if $R=\{0\}$. Retain the rows whose slack stays constant along $r_0$, that is, $J=\{i:a_i^\top r_0=0\}$.

For each $i\in J$, we have $0=a_i^\top r_0=\sum_{e\in E}a_i^\top e$. Since $a_i^\top e\le0$ for every $e\in E$, it follows that $a_i^\top e=0$ for every $e\in E$. As $E$ generates $R$, we have $a_i^\top r=0$ for every $r\in R$. By the Minkowski--Weyl theorem, $K=\conv(V)+R$ for some finite nonempty set $V$. Thus $s_i(v+r)=s_i(v)$ for $v\in\conv(V)$ and $r\in R$. Since $s_i$ is affine, $B_i=\max_{v\in V}s_i(v)<\infty$.

Let $z\in\Z^n$ satisfy $a_i^\top z\le b_i$ for every $i\in J$. For $i\in J$, we have $a_i^\top r_0=0$, so $a_i^\top(z+tr_0)=a_i^\top z\le b_i$ for every $t\ge0$. For $i\notin J$, we have $a_i^\top r_0<0$, so $a_i^\top(z+tr_0)\le b_i$ whenever
\[
t\ge \frac{a_i^\top z-b_i}{-a_i^\top r_0}.
\]
Choose a nonnegative integer $t$ at least as large as each of these finitely many lower bounds. Then $z+tr_0$ satisfies every defining inequality of $P$. Since $z$ and $r_0$ are integral, $z+tr_0\in P\cap\Z^n$, where $g\le0$. As $c^\top r_0=0$, we obtain $g(z)=g(z+tr_0)\le0$.
\end{proof}

Lemma~\ref{lem:rows} shows that the region where $g>0$ and $s_i>0$ for every $i\in J$ contains no integer point. To find a split whose two branches stay away from a point $x$ of positive height, we first bound how far $x$ can move while remaining in this region. Let $h=H/4$ and consider $x\in Q$ with $g(x)\ge h$. The matrix $D_x$ has rows $a_i^\top/s_i(x)$ for $i\in J$ and last row $c^\top/(g(x)+h)$. For $y\in\R^n$,
\begin{equation}\label{eq:localnorm}
\|D_x(y-x)\|_2^2=\sum_{i\in J}\left(\frac{s_i(y)-s_i(x)}{s_i(x)}\right)^2
+\left(\frac{g(y)-g(x)}{g(x)+h}\right)^2.
\end{equation}
If $\|D_x(y-x)\|_2<1/2$, each ratio in~\eqref{eq:localnorm} has absolute value less than $1/2$. Hence $s_i(y)>s_i(x)/2>0$ for every $i\in J$, and $g(y)>(g(x)-h)/2\ge0$.

We now consider the images of the integer points in $W=D_x\R^n$. For $z\in\Z^n$, each $a_i^\top z$ and $c^\top z$ is an integer. Thus the coordinate of $D_xz$ corresponding to row $i\in J$ is an integer multiple of $1/s_i(x)$, and its last coordinate is an integer multiple of $1/(g(x)+h)$. Two distinct points of $\Lambda=D_x\Z^n$ differ in at least one coordinate. Their distance is therefore at least the smallest of these positive spacings. Hence $\Lambda$ is discrete. It is the additive group generated by the columns of $D_x$, since $D_xz$ is their linear combination with coefficients $z_1,\ldots,z_n$. Those same columns span $W$ over $\mathbb R$, so $\Lambda$ is a lattice spanning $W$.

Suppose a point of $\Lambda$ were at distance less than $1/2$ from $D_xx$. That point has the form $D_xz$ for some $z\in\Z^n$, so $\|D_x(z-x)\|_2<1/2$. By~\eqref{eq:localnorm}, $s_i(z)>0$ for every $i\in J$ and $g(z)>0$. Since $z$ is integral, $s_i(z)>0$ implies $a_i^\top z\le b_i$. Thus $z$ satisfies every inequality indexed by $J$, and Lemma~\ref{lem:rows} implies $g(z)\le0$, a contradiction. Hence $\dist(D_xx,\Lambda)\ge1/2$.

The bound $\dist(D_xx,\Lambda)\ge1/2$ measures the distance to the nearest lattice point. To construct a split, we need a direction in which $D_xx$ lies between two consecutive parallel hyperplanes and stays a positive distance from both. We choose the direction from the dual lattice
\[
\Lambda^*=\{u\in W:\langle u,z\rangle\in\Z\text{ for every }z\in\Lambda\}.
\]
For $0\ne u\in\Lambda^*$, every point of $\Lambda$ lies on one of the hyperplanes $\langle u,w\rangle=k$, $k\in\Z$. The distance from $p\in W$ to the hyperplane indexed by $k$ is $\lvert\langle u,p\rangle-k\rvert/\|u\|_2$. Minimizing over $k\in\Z$ yields the distance to the nearest hyperplane. For $d=\dim W$, H\aa stad~\cite{Hastad} showed that one can choose $u\in\Lambda^*\setminus\{0\}$ so that the distance from $p$ to the nearest of these hyperplanes is at least $\dist(p,\Lambda)/(6d^2+1)$. Banaszczyk~\cite[Theorem~2.3]{Ban} improved the dependence of the denominator on $d$ from quadratic to linear. We use the following consequence of his theorem.

\begin{lemma}\label{rank:lem:dual}
Let $\Lambda$ be a lattice spanning a Euclidean space $W$ of dimension $d\ge1$. For every $p\in W\setminus\Lambda$, there is a nonzero $u\in\Lambda^*$ such that
\[
\frac{\dist(\langle u,p\rangle,\Z)}{\|u\|_2}\ge\frac{\dist(p,\Lambda)}{6d}.
\]
\end{lemma}
\begin{proof}
Suppose $d\ge2$. By Banaszczyk~\cite[Theorem~2.3]{Ban}, there is a vector $u\in\Lambda^*$ satisfying
\[
\|u\|_2\,\dist(p,\Lambda)\le d,\qquad \cos(2\pi\langle u,p\rangle)\le\tfrac1{10}.
\]
Let $k$ be an integer nearest to $\langle u,p\rangle$. If $\lvert\langle u,p\rangle-k\rvert\le1/6$, then $\cos(2\pi\langle u,p\rangle)=\cos(2\pi(\langle u,p\rangle-k))\ge1/2$, contrary to the second inequality. Thus $\dist(\langle u,p\rangle,\Z)>1/6$, which also implies $u\ne0$. The first inequality implies $1/\|u\|_2\ge\dist(p,\Lambda)/d$. Multiplying these bounds proves the claim.

For $d=1$, choose a unit coordinate in which $\Lambda=b\Z$, with $b>0$. Then $u=1/b$ belongs to $\Lambda^*$, and
\[
\frac{\dist(\langle u,p\rangle,\Z)}{\|u\|_2}
=b\min_{k\in\Z}|p/b-k|=\min_{k\in\Z}|p-bk|=\dist(p,\Lambda).\qedhere
\]
\end{proof}

\begin{figure}[tb]
\centering
\begin{tikzpicture}[x=2.04cm,y=2.04cm,>=Stealth,font=\small]
  \fill[gray!11] (-1.15,-.21) rectangle (0,1.21);
  \fill[gray!11] (1,-.21) rectangle (1.33,1.21);
  \draw[inkblue] (0,-.21)--(0,1.21) node[above] {$z_1=0$};
  \draw[inkblue] (1,-.21)--(1,1.21) node[above] {$z_1=1$};
  \draw[cutteal,thick,dashed] (.25,.5) circle[radius=1.02cm];
  \draw[inkblue,densely dotted] (.25,.5)--(0,0);
  \draw[cutteal,very thick] (0,.5)--(.25,.5);
  \foreach \a in {-1,0,1}{\foreach \b in {0,1}{\fill[inkblue] (\a,\b) circle(1.45pt);}}
  \fill[cutteal] (.25,.5) circle(1.9pt) node[right] {$p$};
\end{tikzpicture}
\caption{From $p=(\tfrac14,\tfrac12)$, the nearest point of $\mathbb Z^2$ is at distance $\sqrt5/4$ (dotted), and the nearer branch of $(e_1,0)$ is at distance $\tfrac14$ (solid). The ball of radius $\tfrac12$ about $p$ contains no lattice point but meets the left branch.}\label{fig:lattice}
\end{figure}
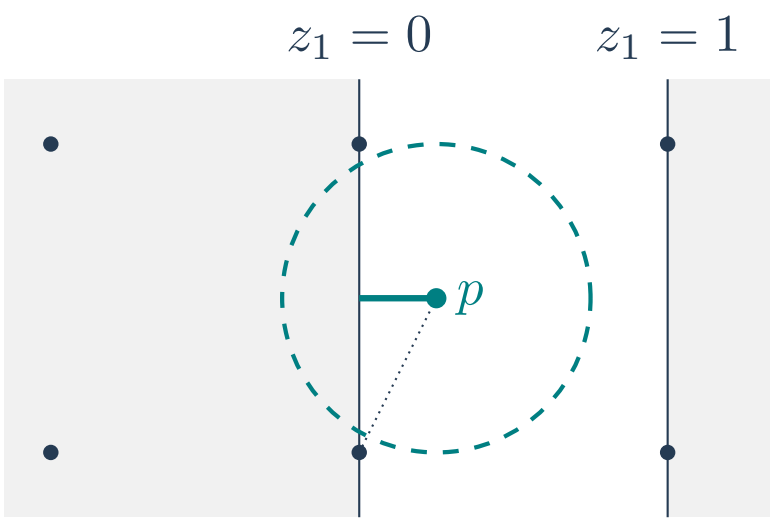

Lemma~\ref{rank:lem:dual} determines a split in the original variables. Take $p=D_xx$ and let $u\in\Lambda^*$ satisfy the lemma. The vector $\pi=D_x^\top u$ is integral because each of its entries is the inner product of $u$ with a column of $D_x$, and each column belongs to $\Lambda$. Since $\dist(p,\Lambda)\ge1/2$, Lemma~\ref{rank:lem:dual} implies that $\pi^\top x=\langle u,p\rangle$ is not an integer. In particular, $\pi\ne0$. Let $k=\lfloor\pi^\top x\rfloor$. For $y$ on either branch of $(\pi,k)$, the interval between $\pi^\top x$ and $\pi^\top y$ contains $k$ or $k+1$. Thus $\lvert\pi^\top(y-x)\rvert\ge\dist(\pi^\top x,\Z)$. By Cauchy--Schwarz and Lemma~\ref{rank:lem:dual},
\begin{equation}\label{rank:eq:sep}
 \|D_x(y-x)\|_2\ge\frac{|\pi^\top(y-x)|}{\|u\|_2}
 \ge\frac{\dist(\pi^\top x,\Z)}{\|u\|_2}\ge\frac1{12n}=: \eta,
\end{equation}
where the last inequality uses $\dist(p,\Lambda)\ge1/2$ and $\dim W\le n$. Figure~\ref{fig:lattice} illustrates the distances from $p$ to the lattice and to the nearer split branch.

To prove that a cut holds at a branch point $y$ with $g(y)<0$, we will use the point where the segment $[x,y]$ crosses $g=0$. This point can lie between the branches, where~\eqref{rank:eq:sep} does not apply. On the entire section $g=0$, however, the height term in~\eqref{eq:localnorm} already suffices. For every $w$ with $g(w)=0$,
\begin{equation}\label{eq:zero}
 \|D_x(w-x)\|_2\ge\frac{|g(w)-g(x)|}{g(x)+h}=\frac{g(x)}{g(x)+h}\ge\tfrac12,
\end{equation}
since $g(x)\ge h$.

\subsection{The halving bound}\label{sec:halving}
We use the distance estimates~\eqref{rank:eq:sep} and~\eqref{eq:zero} to bound the gap between a convex function and its affine tangent at $x$. This gap will allow us to raise the tangent while keeping it below the function on both split branches at nonnegative height and on the entire section $g=0$. For $y\in P$ with $0\le g(y)\le H$, let
\begin{equation}\label{eq:barrier}
t(y)=\frac{g(y)+h}{h},\qquad F(y)=-\sum_{i\in J}\log s_i(y)-\log t(y),
\end{equation}
where $\log$ is the natural logarithm. Since $h=H/4$, we have $1\le t(y)\le5$. Also, $y\in K$, so $1\le s_i(y)\le B_i$ for every $i\in J$. Consequently,
\begin{equation}\label{rank:eq:osc}
-\Delta\le F(y)\le0,\qquad \Delta=\log 5+\sum_{i\in J}\log B_i.
\end{equation}
The same $\Delta$ applies at every halving because $K$, $J$, and the $B_i$ depend only on $P$ and $g$.

Let $T_x(y)=F(x)+\nabla F(x)^\top(y-x)$ be the affine tangent to $F$ at $x$. For each $s$ among the shifted slacks $(s_i)_{i\in J}$ and $t$, let $v_s=(s(y)-s(x))/s(x)$. Since $s(y)=s(x)(1+v_s)$, the change in $-\log s$ from $x$ to $y$ is $-\log(1+v_s)$. The corresponding change in its tangent is $-v_s$, since the derivative of $-\log$ at $s(x)$ is $-1/s(x)$. Summing over $s$, we obtain
\[
F(y)-T_x(y)=\sum_s\bigl(v_s-\log(1+v_s)\bigr).
\]
These are the same relative changes that appear in~\eqref{eq:localnorm}: in particular, $v_t=(g(y)-g(x))/(g(x)+h)$. Thus $\sum_s v_s^2=\|D_x(y-x)\|_2^2$, the quadratic form of the Hessian of $F$ at $x$ evaluated on $y-x$~\cite[Section~11.2]{BV04}. To bound the tangent gap from below, it remains to compare $v-\log(1+v)$ with $v^2$.

Let $\psi(v)=v-\log(1+v)$ for $v>-1$. Since $\psi(0)=\psi'(0)=0$ and $\psi''(v)=1/(1+v)^2>0$, the minimum of $\psi$ is zero. On $(-1,1]$, we have $\psi''(v)\ge1/4$, so Taylor's theorem at zero implies $\psi(v)\ge v^2/8$. For $v>1$, we have $\psi'(v)=v/(1+v)>0$, so $\psi(v)\ge\psi(1)=1-\log 2>1/8$.

Let $y$ lie on either branch of the split $(\pi,k)$, with $0\le g(y)\le H$. By~\eqref{rank:eq:sep}, $\|D_x(y-x)\|_2\ge\eta$. Squaring and using~\eqref{eq:localnorm}, we obtain
\[
\sum_s v_s^2=\|D_x(y-x)\|_2^2\ge\eta^2.
\]
If every $v_s\le1$, then $\psi(v_s)\ge v_s^2/8$ for each $s$. Summing these inequalities yields
\[
F(y)-T_x(y)=\sum_s\psi(v_s)\ge\frac18\sum_s v_s^2\ge\frac{\eta^2}{8}.
\]
If some $v_s>1$, then $\psi(v_s)>1/8$. All other terms are nonnegative, so $F(y)-T_x(y)>1/8\ge\eta^2/8$, since $\eta=1/(12n)\le1$. In both cases the only property of the split branch used was the distance bound in~\eqref{rank:eq:sep}. Thus, for every $y\in P$ with $0\le g(y)\le H$,
\begin{equation}\label{rank:eq:gain}
\|D_x(y-x)\|_2\ge\eta\quad\Longrightarrow\quad F(y)-T_x(y)\ge4\gamma,\qquad \gamma=\frac{\eta^2}{32}.
\end{equation}
In particular,~\eqref{rank:eq:gain} also applies on the entire section $g=0$, since~\eqref{eq:zero} bounds the distance there by at least $1/2$. We use this gap in successive split rounds to lower the upper bound on $F$ at points with $g\ge H/2$. Once that upper bound falls below $-\Delta$,~\eqref{rank:eq:osc} excludes every such point.

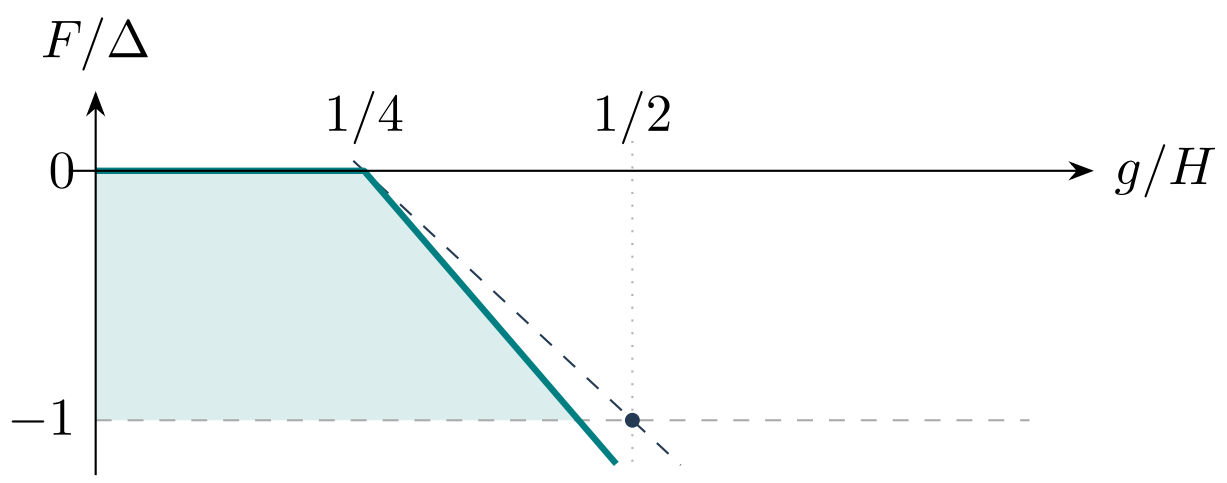
\begin{figure}[tb]
\centering
\begin{tikzpicture}[x=7.1cm,y=1.65cm,font=\small,line join=round,>=Stealth]
 \fill[cutteal!14] (0,-1)--(0,0)--(.25,0)--(.45,-1)--cycle;
 \draw[gray!65,dashed] (0,-1)--(.87,-1);
 \draw[gray!60,dotted] (.5,.12)--(.5,-1.18);
 \draw[inkblue,dashed] (.24,.04)--(.545,-1.18);
 \draw[cutteal,very thick] (0,0)--(.25,0)--(.485,-1.175);
 \draw[->] (-.025,0)--(.93,0) node[right] {$g/H$};
 \draw[->] (0,-1.22)--(0,.32) node[above] {$F/\Delta$};
 \node[left] at (0,-1) {$-1$};
 \node[left] at (0,0) {$0$};
 \node[above] at (.25,.025) {$1/4$};
 \node[above] at (.5,.025) {$1/2$};
 \fill[inkblue] (.5,-1) circle(1.4pt);
\end{tikzpicture}
\caption{In the normalized coordinates $(g/H,F/\Delta)$, the dashed line of slope $-4$ passes through $(\tfrac12,-1)$. When $j\gamma/\Delta>4$, the upper bound in~\eqref{rank:eq:ind} is steeper, so the shaded region satisfying~\eqref{rank:eq:osc} and~\eqref{rank:eq:ind} lies entirely in $g<H/2$.}\label{fig:halving}
\end{figure}

\begin{lemma}\label{lem:halving}
If $g\le H$ on $Q$, then $g\le H/2$ on $\SC^b(Q)$ for $b=1+\lfloor4\Delta/\gamma\rfloor$.
\end{lemma}
\begin{proof}
Let $Q_j=\SC^j(Q)$. We seek an upper bound on $F$ that decreases with each round at points of height at least $H/2$. The split in~\eqref{rank:eq:sep} is available when $g\ge h$, where $h=H/4$. To account for this threshold, let $\phi(y)=(g(y)-h)/H$. Then $\phi(y)\le0$ when $g(y)\le h$, and $\phi(y)\ge1/4$ when $g(y)\ge H/2$. We prove by induction on $j$ that
\begin{equation}\label{rank:eq:ind}
 F(y)+j\gamma\phi(y)\le0\qquad\text{for every }y\in Q_j\text{ with }g(y)\ge0.
\end{equation}
For $j=0$, this is the upper bound $F(y)\le0$ in~\eqref{rank:eq:osc}. When $0\le g(y)\le h$, both terms on the left of~\eqref{rank:eq:ind} are nonpositive for every $j$. When $g(y)\ge H/2$,~\eqref{rank:eq:ind} requires $F(y)\le-j\gamma/4$.

Assume~\eqref{rank:eq:ind} holds on $Q_j$, and let $x\in Q_j$ with $g(x)\ge0$ satisfy
\[
F(x)+(j+1)\gamma\phi(x)>0.
\]
Since $F(x)\le0$ by~\eqref{rank:eq:osc}, we must have $\phi(x)>0$. Thus $g(x)>h$ by the definition of $\phi$, and the split $(\pi,k)$ in~\eqref{rank:eq:sep} applies at $x$. Since $\phi$ is affine, the tangent at $x$ to $F+j\gamma\phi$ is $T_x+j\gamma\phi$.

Let $y\in Q_j$ lie on either branch of $(\pi,k)$, with $g(y)\ge0$. Since $Q_j\subseteq Q$, we have $g(y)\le H$. By~\eqref{rank:eq:sep}, $\|D_x(y-x)\|_2\ge\eta$, so~\eqref{rank:eq:gain} implies $F(y)-T_x(y)\ge4\gamma$. Adding $j\gamma\phi(y)$ to both sides of $T_x(y)\le F(y)-4\gamma$ yields
\[
T_x(y)+j\gamma\phi(y)\le F(y)+j\gamma\phi(y)-4\gamma\le-4\gamma,
\]
where the last inequality follows from the induction hypothesis~\eqref{rank:eq:ind}. Thus we can raise the tangent by $2\gamma$ and still have $\ell_x(y)\le-2\gamma$ on both branches at nonnegative height, where
\begin{equation}\label{eq:cut}
 \ell_x(y)=T_x(y)+j\gamma\phi(y)+2\gamma.
\end{equation}
The coefficients of $\ell_x$ may be irrational. Validity on both branches suffices here because every affine inequality valid on both branches holds on their closed convex hull. Section~\ref{sec:encoding} constructs a rational cut for the same split that also excludes $x$.

At $x$, we have $T_x(x)=F(x)$ and $\phi(x)\le3/4$, since $g(x)\le H$. Hence
\begin{equation}\label{eq:positive-margin}
 \ell_x(x)=\underbrace{F(x)+(j+1)\gamma\phi(x)}_{>0}+(2-\phi(x))\gamma>\gamma.
\end{equation}
Thus $x$ violates $\ell_x\le0$, whereas every $y\in Q_j$ on either branch with $g(y)\ge0$ satisfies
\begin{equation}\label{eq:negative-margin}
 \ell_x(y)\le F(y)+j\gamma\phi(y)-2\gamma\le-2\gamma.
\end{equation}

The bound $\ell_x(w)\le-2\gamma$ also holds for every $w\in P$ with $g(w)=0$, including points between the two branches. Indeed,~\eqref{eq:zero} implies $\|D_x(w-x)\|_2\ge1/2\ge\eta$, so $T_x(w)\le F(w)-4\gamma$ by~\eqref{rank:eq:gain}. Since $F(w)\le0$ by~\eqref{rank:eq:osc} and $\phi(w)=-1/4$, equation~\eqref{eq:cut} implies
\[
\ell_x(w)\le F(w)-4\gamma-\frac{j\gamma}{4}+2\gamma\le-2\gamma.
\]

Now let $y\in Q_j$ lie on either branch with $g(y)<0$. Since $g(x)>h>0$, the segment $[x,y]$ contains a point $w$ with $g(w)=0$. If $\ell_x(y)>0$, then $\ell_x(x)>0$ by~\eqref{eq:positive-margin}, and affineness implies $\ell_x(w)>0$, contradicting $\ell_x(w)\le-2\gamma$. Thus $\ell_x(y)\le0$ on both branches at negative height as well. Together with~\eqref{eq:negative-margin}, this proves that $\ell_x\le0$ is a split cut for $Q_j$. For every $x\in Q_j$ with $g(x)\ge0$ and
\[
F(x)+(j+1)\gamma\phi(x)>0,
\]
we have constructed a split cut for $Q_j$ that excludes $x$. Since $Q_{j+1}=\SC(Q_j)$ satisfies every split cut for $Q_j$, every $y\in Q_{j+1}$ with $g(y)\ge0$ satisfies
\[
F(y)+(j+1)\gamma\phi(y)\le0.
\]
This proves the induction step.

Suppose $y\in Q_j$ satisfies $g(y)\ge H/2$. Since $h=H/4$, we have $\phi(y)=(g(y)-h)/H\ge1/4$. The lower bound in~\eqref{rank:eq:osc} and the upper bound in~\eqref{rank:eq:ind} therefore imply
\[
-\Delta\le F(y)\le-j\gamma\phi(y)\le-\frac{j\gamma}{4}.
\]
Such a point can exist only if $j\le4\Delta/\gamma$. Since $b=1+\lfloor4\Delta/\gamma\rfloor>4\Delta/\gamma$, every $y\in Q_b=\SC^b(Q)$ satisfies $g(y)<H/2$, as illustrated in Figure~\ref{fig:halving}.
\end{proof}

\section{The rank bound}\label{sec:rankbound}
Lemma~\ref{lem:halving} bounds the number of split rounds needed to halve an upper bound on $g$. To prove Theorem~\ref{thm:generalrank} for $P_I\ne\varnothing$, we need two estimates for each defining inequality $g\le0$ of $P_I$ that is violated on $P$. The initial height $H_0=\lceil\max_P g\rceil$ determines how many halvings are needed to reach $g<1$. The quantity $\Delta$ in~\eqref{rank:eq:osc} determines how many rounds each halving takes. We bound both in terms of the input length $L$. Let $N=L+1$. The dense encoding counts every coefficient and right-hand side, including zeros, and each entry in the integer description of $P$ requires $O(N)$ bits. If $m=0$, then $P=\mathbb R^n=\conv(\mathbb Z^n)$. If $n=0$, then $P$ is either empty or consists of the unique point of $\mathbb R^0$, which is integral. In either case $P=P_I$, so the split rank is zero. We may therefore assume $m,n\ge1$. The dense encoding records $n$ coefficients and one right-hand side for each of the $m$ inequalities, so $L\ge m(n+1)$. In particular, $m,n\le N$.

\begin{lemma}\label{lem:sizes}
Suppose $P_I\ne\varnothing$. Then $P_I$ has a finite description by inequalities $g\le0$ with $O(N^3)$ bits per coefficient, including the constant term. For each inequality in this description that is violated on $P$, the maximum of $g$ on $P$ is finite. With $H_0=\lceil\max_P g\rceil$ and $B_i=\max_{y\in K}s_i(y)$ for the rows $i\in J$ of Lemma~\ref{lem:rows}, we have
\[
 \log_2(1+H_0)=O(N^3),\qquad \log B_i=O(N^4)\quad(i\in J).
\]
If $P\ne\varnothing$ and $P_I=\varnothing$, the same bounds hold for $g\equiv1$ and $H_0=1$.
\end{lemma}

\begin{proof}
We use a determinant bound three times. First, we bound the coordinates of the point and recession generators of $P$. These bounds determine a box $B$ large enough that $P_I=\conv(P\cap B\cap\Z^n)+\rec(P)$. Second, we use this description to express validity of an inequality $g\le0$ as a system of linear inequalities in the coefficients of $g$. Applying the determinant bound to this system bounds the coefficients of a defining description of $P_I$. We then bound $H_0$ by evaluating each defining function $g$ on the point generators of $P$. Third, we apply the determinant bound to $K=P\cap\{g\ge0\}$. By Lemma~\ref{lem:rows}, each retained slack $s_i$ is constant along every recession direction of $K$, so its maximum $B_i$ is attained at a point generator of $K$.

Let $C=\{w\in\mathbb R^d:Aw\le0\}$, where $A$ is an integer matrix and $|A_{ij}|\le2^q$ for some $q\ge1$. We show that $C$ has a finite set of integer generators, each of whose coordinates requires at most
\begin{equation}\label{eq:det}
O\bigl(d(q+\log(d+1))\bigr)\quad\text{bits}.
\end{equation}
The cone $C$ may contain a line. To obtain cones generated by extreme rays, intersect $C$ with each coordinate orthant. In particular, for each $\sigma\in\{-1,1\}^d$, let
\[
C_\sigma=\{w\in C:\sigma_j w_j\ge0\text{ for }j=1,\ldots,d\}.
\]
These cones cover $C$, since every vector has a choice of coordinate signs. Each $C_\sigma$ is pointed: if $w$ and $-w$ both belong to $C_\sigma$, then $\sigma_jw_j=0$ for every $j$, so $w=0$. Let $r\ne0$ lie on an extreme ray of $C_\sigma$. Form $M$ from the rows $A_i$ for which $A_i r=0$ and the rows $\sigma_j e_j^\top$ for which $r_j=0$. These are precisely the constraints defining $C_\sigma$ that hold with equality at $r$. Since $Mr=0$ and $r\ne0$, we have $\operatorname{rank}(M)\le d-1$. If the rank were smaller, the nullspace of $M$ would contain a vector $v$ linearly independent of $r$. Then $M(r\pm\varepsilon v)=0$, and for sufficiently small $\varepsilon>0$, all remaining inequalities would still hold because they are strict at $r$. Thus $r$ would be the midpoint of two points of $C_\sigma$ outside its ray, contradicting extremality. Hence $\operatorname{rank}(M)=d-1$.

For $d\ge2$, choose $d-1$ linearly independent rows of $M$ and let $B$ be the matrix with these rows. Since $Br=0$ and $B$ has rank $d-1$, its nullspace is the line spanned by $r$. To find an integer vector on this line, let $B_{\widehat j}$ be the square matrix obtained by deleting column $j$ of $B$, and define
\[
u_j=(-1)^{j+1}\det B_{\widehat j},\qquad j=1,\ldots,d.
\]
To verify $Bu=0$, consider any row $b^\top$ of $B$. Placing $b^\top$ above $B$ produces a square matrix with two identical rows, so its determinant is zero. Expanding along the first row yields $\sum_{j=1}^d b_j(-1)^{j+1}\det B_{\widehat j}=b^\top u$. Thus $b^\top u=0$ for every row of $B$. Since $B$ has rank $d-1$, some $d-1$ columns are linearly independent, so at least one determinant $\det B_{\widehat j}$ is nonzero. Thus $u$ is a nonzero integer vector in $\linspan(r)$. Since $r\ne0$, either $u$ or $-u$ generates the ray containing $r$.

Each row of $B$ is either a row of $A$ or a coordinate row $\sigma_j e_j^\top$, so every entry has absolute value at most $2^q$. Each of the $d-1$ rows of $B_{\widehat j}$ therefore has Euclidean norm at most $\sqrt d\,2^q$. Hadamard's inequality bounds the absolute determinant by the product of these row norms. Hence
\[
|u_j|\le (\sqrt d\,2^q)^{d-1}.
\]
The binary length of $u_j$, including its sign, is consequently at most $(d-1)(q+\tfrac12\log_2 d)+O(1)$, as required in~\eqref{eq:det}. For $d=1$, the positive and negative rays have generators $1$ and $-1$, respectively, so~\eqref{eq:det} also holds.

Finally, each $C_\sigma$ is a pointed polyhedral cone, so it is generated by its finitely many extreme rays. Every $w\in C$ belongs to some $C_\sigma$ and is therefore a nonnegative combination of the generators chosen for that cone. Conversely, every chosen generator lies in $C$, and $C$ contains all their nonnegative combinations. Taking the generators over the finitely many sign vectors $\sigma$ therefore gives a finite generating set for $C$ satisfying~\eqref{eq:det}.

By the Minkowski--Weyl theorem, $P=\{y:Ay\le b\}$ has finite sets of point generators $V$ and recession generators $E$ such that
\begin{equation}\label{eq:gen}
 P=\conv(V)+\cone(E),\qquad \rec(P)=\cone(E).
\end{equation}
To bound their coordinates, apply~\eqref{eq:det} to the homogenized cone $\widehat P=\{(y,t):Ay\le tb,\ t\ge0\}$. From its integer generators, take $v/t\in V$ for each $(v,t)$ with $t>0$, and $e\in E$ for each $(e,0)$. The system defining $\widehat P$ has $n+1$ variables and $O(N)$ bits per coefficient. Since $n\le N$,~\eqref{eq:det} bounds each generator coordinate by $O(N^2)$ bits. As $t$ is a positive integer for each point generator, $\|v/t\|_\infty\le\|v\|_\infty$. Thus every vector in $V\cup E$ has infinity norm at most $2^{O(N^2)}$.

Suppose $P_I\ne\varnothing$. Let $B=1+\max_{v\in V}\|v\|_\infty+n\max_{e\in E}\|e\|_\infty$, taking the second maximum to be zero if $E=\varnothing$. For $Z_B=P\cap\mathbb Z^n\cap[-B,B]^n$, we claim that
\begin{equation}\label{eq:box}
P_I=\conv(Z_B)+\rec(P).
\end{equation}
Let $y\in P\cap\mathbb Z^n$. By~\eqref{eq:gen} and conic Carath\'eodory~\cite[Corollary~7.1i]{Schrijver86}, $y=v+\sum_{j=1}^k\lambda_j e_j$ for some $v\in\conv(V)$, $e_j\in E$, and $\lambda_j\ge0$, with $k\le n$. The point
\[
z=y-\sum_{j=1}^k\lfloor\lambda_j\rfloor e_j
=v+\sum_{j=1}^k(\lambda_j-\lfloor\lambda_j\rfloor)e_j
\]
is integral by the first equality and belongs to $P$ by the second. Each remaining coefficient is less than one, so $\|z\|_\infty<B$ by the choice of $B$. Thus $z\in Z_B$ and $y-z\in\rec(P)$. Taking convex hulls proves the forward inclusion in~\eqref{eq:box}.

Conversely, $\conv(Z_B)\subseteq P_I$, and $\rec(P)=\rec(P_I)$ because $P_I\ne\varnothing$~\cite{Schrijver86}. Hence $\conv(Z_B)+\rec(P)\subseteq P_I$, proving~\eqref{eq:box}.

By~\eqref{eq:gen} and~\eqref{eq:box}, an inequality $\alpha^\top y\le\beta$ is valid for $P_I$ if and only if
\[
\alpha^\top z\le\beta\quad(z\in Z_B),\qquad \alpha^\top e\le0\quad(e\in E).
\]
These inequalities define a cone in the $n+1$ variables $(\alpha,\beta)$. Each coordinate of $z\in Z_B$ is an integer of absolute value at most $B$. The coordinate bound following~\eqref{eq:gen} and the definition of $B$ imply that $B\le2^{O(N^2)}$ and $\|e\|_\infty\le2^{O(N^2)}$ for every $e\in E$. Thus each coefficient of this system has $O(N^2)$ bits. By~\eqref{eq:det}, this cone has a finite set of integer generators with $O(N^3)$ bits per coordinate. The corresponding inequalities describe $P_I$. Indeed, by separation, every point outside $P_I$ violates a valid inequality. Since that inequality is a nonnegative combination of the generator inequalities, the point must violate at least one of them.

Let $g(y)=c^\top y+c_0$, where $g\le0$ is one of these defining inequalities and $g$ is positive somewhere on $P$. Validity on $P_I$ implies $c^\top r\le0$ for every $r\in\rec(P_I)=\rec(P)$. Thus~\eqref{eq:gen} implies $\max_P g=\max_{v\in V}g(v)$. The coefficients of $g$ have $O(N^3)$ bits, and $\|v\|_\infty\le2^{O(N^2)}$ for every $v\in V$. Since $n\le N$, it follows that $\log_2(1+H_0)=O(N^3)$. To bound the retained slacks, apply the same construction to $K=P\cap\{g\ge0\}$. Its homogenized cone has $n+1$ variables, and the additional inequality $g\ge0$ has $O(N^3)$ bits per coefficient. By~\eqref{eq:det}, $K$ therefore has
\begin{equation}\label{eq:Ksize}
\text{point generators with }O(N^4)\text{ bits per coordinate}.
\end{equation}
For each $i\in J$, Lemma~\ref{lem:rows} shows that $B_i=\max_K s_i$ is attained at one of these point generators. Since $s_i(y)=b_i+1-a_i^\top y$ has $O(N)$ bits per coefficient,~\eqref{eq:Ksize} implies $\log B_i=O(N^4)$. If $P\ne\varnothing$ and $P_I=\varnothing$, take $g\equiv1$. Then $H_0=1$ and $K=P$, so the bound on $B_i$ follows directly from the point-generator bound following~\eqref{eq:gen}.
\end{proof}

\begin{proof}[Proof of Theorem~\ref{thm:generalrank}]
Suppose $P_I\ne\varnothing$, and let $g(y)=c^\top y+c_0$, where $g\le0$ is a defining inequality from Lemma~\ref{lem:sizes} that is violated on $P$. Since $g\le0$ on $P_I$ and $g>0$ somewhere on $P$, we have $c\ne0$. The split $(c,-c_0)$ has branches $g\le0$ and $g\ge1$. To enforce $g\le0$ in one split round, it therefore suffices to reach a relaxation on which $g<1$.

Let $H_0=\lceil\max_P g\rceil$ and $H_\nu=H_0/2^\nu$. We claim that $g\le H_\nu$ on $\SC^{\nu b}(P)$ for every $\nu\ge0$. At $\nu=0$, this follows from the definition of $H_0$. Suppose the claim holds at $\nu$. Apply Lemma~\ref{lem:halving} with $Q=\SC^{\nu b}(P)$ and $H=H_\nu$. Then $g\le H_\nu/2=H_{\nu+1}$ on $\SC^b(Q)=\SC^{(\nu+1)b}(P)$. The number $b=1+\lfloor4\Delta/\gamma\rfloor$ is the same at each application: the retained rows $J$ and their slack maxima $B_i$ in~\eqref{rank:eq:osc} are determined by $P$ and $g$, and $\gamma$ in~\eqref{rank:eq:gain} depends only on $n$. This proves the claim. For $a=1+\lfloor\log_2 H_0\rfloor$, we have $2^a>H_0$. Thus $g\le H_a<1$ on $\SC^{ab}(P)$, and the split $(c,-c_0)$ enforces $g\le0$ in round $ab+1$.

We now bound $ab+1$ in terms of the input length, uniformly over the defining inequalities of $P_I$. Lemma~\ref{lem:sizes} implies $\log_2(1+H_0)=O(N^3)$, so $a=O(N^3)$. To bound $b$, recall from~\eqref{rank:eq:osc} that $\Delta=\log 5+\sum_{i\in J}\log B_i$. There are at most $m$ terms in the sum, each $O(N^4)$ by Lemma~\ref{lem:sizes}, so $\Delta=O(mN^4)$. By~\eqref{rank:eq:sep} and~\eqref{rank:eq:gain}, $\gamma^{-1}=32(12n)^2$. Thus $b=1+\lfloor4\Delta/\gamma\rfloor=O(mn^2N^4)$. The dense input records all $mn$ matrix entries, so $mn\le L$. Since $n\le N$, we have $mn^2\le Ln\le N^2$. It follows that each halving uses $O(N^6)$ rounds and that $ab+1=O(N^9)$. The constants are uniform over the defining rows. Hence there is an integer $R=O(N^9)$ such that every row holds on $\SC^R(P)$. Indeed, rows valid on $P$ remain valid throughout, and a row derived at round $r\le R$ holds on $\SC^R(P)\subseteq\SC^r(P)$. Since these inequalities define $P_I$, we have $\SC^{R}(P)\subseteq P_I$. Conversely, every split cut preserves all integer points of $P$ and therefore their convex hull $P_I$. Thus $\SC^{R}(P)=P_I$.

Suppose $P\ne\varnothing$ and $P_I=\varnothing$. Since $P$ contains no integer point, the inequality $g\le0$ with $g\equiv1$ is valid at every integer point of $P$. Apply Lemma~\ref{lem:halving} with $Q=P$ and $H=1$. After $b$ rounds, every point of $\SC^b(P)$ must satisfy $g\le1/2$. Since $g\equiv1$, this closure is empty. Lemma~\ref{lem:sizes} bounds $\log B_i$ by $O(N^4)$ also in this case, so the same estimate $b=O(N^6)$ applies. If $P=\varnothing$, then $P=P_I$ before any round.
\end{proof}

\section{Encoding the individual cuts}\label{sec:encoding}
For a single split closure, each facet admits a defining inequality of encoding length polynomial in the input size~\cite[p.~347]{BCDJplane} (see also~\cite[Theorem~5.5 and Lemma~5.9]{CCZ14}). We show that one polynomial in the original input length bounds the size of both the cuts and their justifying splits throughout the rounds of Theorem~\ref{thm:generalrank}. In the proof of Lemma~\ref{lem:halving}, the coefficients of $\ell_x$ in~\eqref{eq:cut} depend on $F(x)$ and $\nabla F(x)$ and can be irrational. The strict bounds in~\eqref{eq:positive-margin} and~\eqref{eq:negative-margin} allow us to round these coefficients: at $x$, the cut value exceeds $\gamma$, while on either branch at nonnegative height it is at most $-2\gamma$. We will use these bounds to construct a rational inequality that still excludes $x$ and is valid on both branches of the same split. The precision used in this rounding must therefore be uniform over the points $x$, defining rows, and rounds.

\begin{proposition}\label{prop:encoding}
There is a polynomial $p$ such that $O((L+1)^9)$ rounds suffice to reach $P_I$ when each round imposes all rational cuts that are valid on both branches of a split in the current relaxation and for which the cut and its justifying split each have encoding length at most $p(L)$.
\end{proposition}
\begin{proof}
Let $x$ be a point in the current relaxation with $g(x)\ge0$ that violates~\eqref{rank:eq:ind} at index $j+1$, so that
\[
F(x)+(j+1)\gamma\phi(x)>0.
\]
We first bound the encoding length of the normal $\pi=D_x^\top u$ of the split chosen at $x$. By~\eqref{rank:eq:sep},
\[
\eta\|u\|_2\le\dist(\pi^\top x,\Z)\le\tfrac12.
\]
Since $\eta=1/(12n)$, we have $\|u\|_2\le6n$. For the $l$th coordinate vector $e_l$, Cauchy--Schwarz implies
\[
|\pi_l|=|\langle u,D_xe_l\rangle|\le6n\|D_xe_l\|_2.
\]

Let $M$ be the largest absolute entry among the retained row normals and $c$. The column $D_xe_l$ has one entry $a_{il}/s_i(x)$ for each $i\in J$ and one entry $c_l/(g(x)+h)$, from the height term, for a total of $|J|+1$ entries. Each $s_i(x)\ge1$. Each halving starts with $H\ge1$, so $g(x)+h\ge H/4\ge1/4$. Thus every entry has absolute value at most $4M$, and
\[
\|D_xe_l\|_2^2=\sum_{i\in J}\left(\frac{a_{il}}{s_i(x)}\right)^2
+\left(\frac{c_l}{g(x)+h}\right)^2\le(|J|+1)(4M)^2.
\]
Since $|J|\le m$, it follows that $|\pi_l|\le24nM\sqrt{m+1}$. Each $a_{il}$ has $O(N)$ bits by the input bound, and each $c_l$ has $O(N^3)$ bits by Lemma~\ref{lem:sizes}. Hence $\log_2(1+M)=O(N^3)$. As $m,n\le N$, each integer $\pi_l$ requires $O(N^3)$ bits, and the $n$ coordinates of $\pi$ together require $O(N^4)$ bits.

It remains to bound the split offset $k=\lfloor\pi^\top x\rfloor$. Although $K=P\cap\{g\ge0\}$ can be unbounded, $\pi^\top y$ is constant along each recession direction of $K$. Indeed, for every $r\in\rec(K)$, the proof of Lemma~\ref{lem:rows} establishes $a_i^\top r=0$ for all $i\in J$ and $c^\top r=0$. Hence $D_xr=0$ and $\pi^\top r=\langle u,D_xr\rangle=0$. If $V_K$ denotes the set of point generators in~\eqref{eq:Ksize}, then $K=\conv(V_K)+\rec(K)$. Since $x\in K$,
\[
|\pi^\top x|\le\max_{v\in V_K}|\pi^\top v|\le n\|\pi\|_\infty\max_{v\in V_K}\|v\|_\infty.
\]
Each coordinate of $\pi$ has $O(N^3)$ bits, and each coordinate of $v\in V_K$ has $O(N^4)$ bits by~\eqref{eq:Ksize}. Since $|k|\le|\pi^\top x|+1$, the offset requires $O(N^4)$ bits as well.

\begin{figure}[htb]
\centering
\begin{tikzpicture}[font=\small,x=.87cm,y=2.25cm,
  original/.style={inkblue,densely dashed,thick},
  rounded/.style={cutteal,line width=1.1pt}]
\foreach \i/\heading in {0/{Round in $y_1,y_2$},1/{Round in the slack $s$}}{
  \begin{scope}[shift={({6.6*\i},0)}]
    \node at (0,1.03) {\heading};
    \ifnum\i=0
      \fill[cutteal!23] (-1,0)--(1.25,.75)--(2.6,.75)--(2.6,0)--cycle;
    \else
      \fill[cutteal!23] (-2.6,0) rectangle (2.6,1/3);
    \fi
    \draw[gray] (-2.6,.75)--(2.6,.75);
    \draw[gray,-{Stealth[length=3pt]}] (-2.6,0)--(2.85,0) node[right,text=black] {$u$};
    \node[left] at (-2.6,.75) {$\tfrac34$};
    \node[left] at (-2.6,.5) {$\tfrac12$};
    \node[left] at (-2.6,0) {$0$};
    \node at (-2.82,.93) {$v$};
    \draw[original] (-2.6,.5)--(2.6,.5);
    \ifnum\i=0
      \draw[rounded] (-1,0)--(1.25,.75);
    \else
      \draw[rounded] (-2.6,1/3)--(2.6,1/3);
    \fi
    \foreach \n in {-1,0,1,2}{\fill[inkblue] (\n,0) circle (1.1pt);}
    \fill[inkblue] (-2,0) circle (2pt) node[below=2pt,text=black] {$z$};
    \draw[inkblue,fill=white,line width=.9pt] (2,2/3) circle (2.1pt)
      node[left=2pt,text=black] {$x$};
  \end{scope}
}
\end{tikzpicture}
\caption{Rounding the same split cut to the nearest quarter. In coordinates $u=y_2$, $v=2y_2-y_1$, the strip is $0\le v\le3/4$ and the original cut is $\tfrac23v-\tfrac13\le0$ (dashed boundary). Shading shows what each rounded cut retains. On the left, rounding in the original coordinates gives $3v\le u+1$, which loses the integer point $z=(-2,0)$ and retains $x=(2,\tfrac23)$. On the right, with $s=1+v$, the cut is $\tfrac23s-1\le0$. Rounding gives $v\le\tfrac13$, retaining the entire branch $v=0$ and excluding $x$.}\label{fig:slack-rounding}
\end{figure}
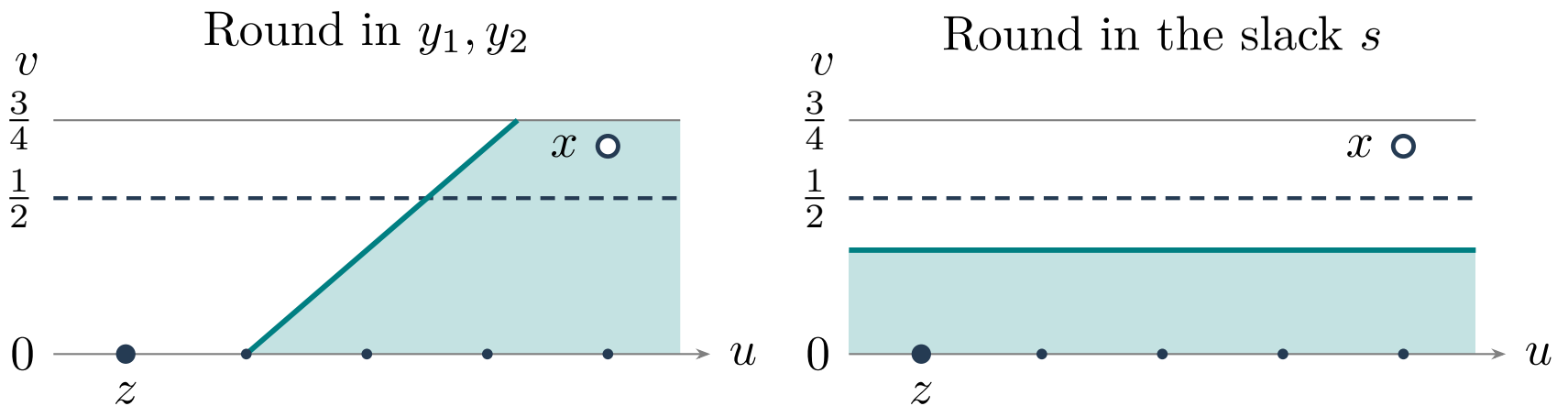

We next round the coefficients of $\ell_x$. Rounding the coefficient of $y_l$ by $\varepsilon$ changes the cut value by $\varepsilon y_l$, for which there need not be a uniform bound on $K$. Figure~\ref{fig:slack-rounding} illustrates how rounding in the original coordinates can exclude an integer point. The slacks in~\eqref{eq:barrier}, however, satisfy $1\le s_i(y)\le B_i$ and $1\le t(y)\le5$ whenever $y\in P$ and $0\le g(y)\le H$. We therefore express $\ell_x$ in these variables. For each slack $s$, its contribution to the tangent $T_x(y)$ is $-\log s(x)+1-s(y)/s(x)$. Summing over the $|J|+1$ slacks and substituting $\phi(y)=t(y)/4-1/2$ into~\eqref{eq:cut}, we obtain
\begin{align}
 \ell_x(y)&=C_x-\sum_{i\in J}\frac{s_i(y)}{s_i(x)}
              +\left(\frac{j\gamma}{4}-\frac{1}{t(x)}\right)t(y),\label{eq:expanded}\\
 C_x&=F(x)+|J|+1+2\gamma-j\gamma/2.\notag
\end{align}
Round the constant and each coefficient of $s_i(y)$ and $t(y)$ in~\eqref{eq:expanded} to a nearest multiple of $2^{-q}$. Let $\delta_0,\delta_i,\delta_t$ denote the rounded value minus the original value for the constant, the coefficient of $s_i(y)$, and the coefficient of $t(y)$, respectively. Consecutive multiples of $2^{-q}$ are $2^{-q}$ apart, so each rounding error has absolute value at most $2^{-q}/2$. The weaker bound $2^{-q}$ will suffice. For the resulting affine function $\widetilde\ell_x$, subtracting~\eqref{eq:expanded} gives
\[
\widetilde\ell_x(y)-\ell_x(y)=\delta_0+\sum_{i\in J}\delta_i s_i(y)+\delta_t t(y).
\]
If $y\in P$ and $0\le g(y)\le H$, then $y\in K$, so $1\le s_i(y)\le B_i$ for each $i\in J$. Also, $t(y)=1+g(y)/h\le1+H/h=5$, since $h=H/4$. By the triangle inequality,
\begin{equation}\label{eq:precision}
\begin{aligned}
|\widetilde\ell_x(y)-\ell_x(y)|&\le|\delta_0|+\sum_{i\in J}|\delta_i|\,s_i(y)+|\delta_t|\,t(y)\\
&\le2^{-q}+\sum_{i\in J}2^{-q}B_i+5\cdot2^{-q}
=2^{-q}\left(6+\sum_{i\in J}B_i\right).
\end{aligned}
\end{equation}

We choose $q$ so that the error in~\eqref{eq:precision} is at most $\gamma/4$. This preserves the strict signs in~\eqref{eq:positive-margin} and~\eqref{eq:negative-margin}: a value greater than $\gamma$ remains greater than $3\gamma/4$, and a value at most $-2\gamma$ remains at most $-7\gamma/4$. The required inequality is
\[
2^{-q}\left(6+\sum_{i\in J}B_i\right)\le\frac{\gamma}{4},
\]
or equivalently $2^q\ge4(6+\sum_{i\in J}B_i)/\gamma$. Taking the base-two logarithm and rounding up, we may therefore choose
\[
q=\left\lceil\log_2\!\left(\frac{4(6+\sum_{i\in J}B_i)}{\gamma}\right)\right\rceil.
\]

We now verify that $\widetilde\ell_x\le0$ is valid on both branches of $(\pi,k)$ and excludes $x$. By~\eqref{eq:positive-margin} and~\eqref{eq:precision}, $\widetilde\ell_x(x)>3\gamma/4$. By~\eqref{eq:negative-margin} and~\eqref{eq:precision}, $\widetilde\ell_x(y)\le-7\gamma/4$ at every branch point $y$ of nonnegative height in the current relaxation. The same bound holds on $P\cap\{g=0\}$, since $\ell_x\le-2\gamma$ throughout that section. Now suppose a branch point $y$ with $g(y)<0$ satisfies $\widetilde\ell_x(y)>0$. Since $g(x)>0$, the segment $[x,y]$ meets $P\cap\{g=0\}$ at some point $w$. The affine function $\widetilde\ell_x$ is positive at both endpoints, so $\widetilde\ell_x(w)>0$, a contradiction. Thus the rounded inequality is a split cut for the current relaxation.

We first bound the number of bits in the rounded coefficients of~\eqref{eq:expanded}. By Lemma~\ref{lem:sizes}, there is an absolute constant $C$ such that $B_i\le2^{CN^4}$ for every $i\in J$. Since $|J|\le N$, we have $\sum_{i\in J}B_i\le N2^{CN^4}$. Substituting this estimate and $\gamma^{-1}=O(N^2)$ into the choice of $q$ following~\eqref{eq:precision}, we obtain
\[
q\le1+\log_2\!\left(\frac4\gamma\right)+\log_2\!\left(6+N2^{CN^4}\right)=O(N^4).
\]

We also need to bound the magnitude of each coefficient before rounding. Since $j\le b-1=\lfloor4\Delta/\gamma\rfloor$, we have $j\gamma\le4\Delta$. The coefficient of $s_i(y)$ in~\eqref{eq:expanded} has absolute value $1/s_i(x)\le1$. The coefficient of $t(y)$ has absolute value at most $j\gamma/4+1/t(x)\le\Delta+1$. For the constant term,~\eqref{rank:eq:osc} implies $|F(x)|\le\Delta$, so
\[
|C_x|\le |F(x)|+|J|+1+2\gamma+j\gamma/2\le3\Delta+|J|+3,
\]
where $\gamma\le1$. Thus $3\Delta+|J|+3$ bounds the magnitude of every coefficient in~\eqref{eq:expanded}. Rounding a coefficient $\alpha$ to a nearest multiple $A/2^q$ means choosing an integer $A$ nearest to $2^q\alpha$. Hence
\[
|A|\le2^q|\alpha|+\tfrac12\le2^q(3\Delta+|J|+3)+1.
\]
Since $\Delta=O(N^5)$, $|J|\le N$, and $q=O(N^4)$, each such integer $A$ requires $O(N^4)$ bits.

It remains to express the rounded cut in the original variables. Denote its integer numerators by $A_0,A_i,A_t$, so that
\[
2^q\widetilde\ell_x(y)=A_0+\sum_{i\in J}A_i s_i(y)+A_t t(y).
\]
At halving $\nu$, we have $h=H_0/2^{\nu+2}$. Substituting $s_i(y)=b_i+1-a_i^\top y$ and $t(y)=1+g(y)/h$, we obtain
\begin{align*}
2^qH_0\,\widetilde\ell_x(y)
={}&H_0\!\left(A_0+A_t+\sum_{i\in J}A_i(b_i+1)\right)+2^{\nu+2}A_t c_0\\
&+\left(2^{\nu+2}A_t c-H_0\sum_{i\in J}A_i a_i\right)^\top y .
\end{align*}
Every quantity on the right has integer coefficients. The factors $H_0$ and $2^{\nu+2}$ each require $O(N^3)$ bits, since Lemma~\ref{lem:sizes} bounds $\log_2(1+H_0)$ by $O(N^3)$ and $\nu\le\lfloor\log_2 H_0\rfloor$. The input coefficients $a_i,b_i$ require $O(N)$ bits per entry, and $c,c_0$ require $O(N^3)$ bits per entry by Lemma~\ref{lem:sizes}. Each product in the display requires $O(N^4)$ bits. Each coefficient is a sum of $O(N)$ such products, so it also requires $O(N^4)$ bits. Since $2^qH_0>0$, the inequality $2^qH_0\,\widetilde\ell_x(y)\le0$ defines the same halfspace as $\widetilde\ell_x(y)\le0$. Its $n+1$ integer coefficients therefore require $O(N^5)$ bits in total.

Choose a polynomial $p(L)$ that bounds the encoding lengths of the rounded cuts and their justifying splits for every defining row $g\le0$, and also for $g\equiv1$ when $P_I=\varnothing$. By Lemma~\ref{lem:sizes}, we may choose the same polynomial to cover the final cut $g\le0$ and its split $(c,-c_0)$ when $c\ne0$. Let $R_0=P$. For each $r\ge0$, obtain $R_{r+1}$ by intersecting $R_r$ with every cut that is valid on both branches of a split in $R_r$ and for which the cut and the split each require at most $p(L)$ bits. There are finitely many such pairs. Since these rounds impose only cuts of bounded length, $R_r$ may contain points outside $\SC^r(P)$. We therefore prove the halving bound directly on $R_r$.

Fix $g$ and suppose a halving starts at $R_s$, with $g\le H$ throughout $R_s$. We prove by induction on $j$, for $0\le j\le b$, that
\[
F(y)+j\gamma\phi(y)\le0\qquad\text{for every }y\in R_{s+j}\text{ with }g(y)\ge0.
\]
At $j=0$, this follows from~\eqref{rank:eq:osc}. Suppose the bound holds at $j<b$, and let $x\in R_{s+j}$ satisfy $g(x)\ge0$ and $F(x)+(j+1)\gamma\phi(x)>0$. Since $F(x)\le0$, we have $\phi(x)>0$, so $g(x)>h$ and the split in~\eqref{rank:eq:sep} applies. For every $y\in R_{s+j}$ on either branch with $g(y)\ge0$, we also have $g(y)\le H$, since $R_{s+j}\subseteq R_s$. Equations~\eqref{rank:eq:sep}, \eqref{rank:eq:gain}, and~\eqref{eq:cut}, together with the induction hypothesis on $R_{s+j}$, imply
\[
\ell_x(y)\le F(y)+j\gamma\phi(y)-2\gamma\le-2\gamma.
\]
Thus~\eqref{eq:precision} ensures that the rounded cut holds at these branch points.

Let $y\in R_{s+j}$ lie on either branch with $g(y)<0$. By~\eqref{eq:positive-margin} and~\eqref{eq:precision}, $\widetilde\ell_x(x)>3\gamma/4$. If $\widetilde\ell_x(y)>0$, then $\widetilde\ell_x$ is positive throughout $[x,y]$ because it is affine. Since $g(x)>0>g(y)$, this segment contains a point $w\in P$ with $g(w)=0$. But~\eqref{eq:zero}, \eqref{rank:eq:gain}, and~\eqref{eq:cut} imply $\ell_x(w)\le-2\gamma$, so~\eqref{eq:precision} implies $\widetilde\ell_x(w)\le-7\gamma/4$, a contradiction. Thus the rounded inequality is valid on both branches of the chosen split in $R_{s+j}$ and excludes $x$. Since both the cut and its split satisfy the length bound $p(L)$, this cut is imposed in forming $R_{s+j+1}$, so $x\notin R_{s+j+1}$. This proves the induction step.

Suppose $y\in R_{s+b}$ satisfies $g(y)\ge H/2$. Then $\phi(y)\ge1/4$, so the induction bound at $j=b$ and the lower bound in~\eqref{rank:eq:osc} imply
\[
-\Delta\le F(y)\le-b\gamma\phi(y)\le-\frac{b\gamma}{4}<-\Delta,
\]
where the strict inequality follows from $b>4\Delta/\gamma$. Thus $g<H/2$ throughout $R_{s+b}$, and we can begin the next halving with upper bound $H/2$. For each defining inequality $g\le0$ violated on $P$, the proof of Theorem~\ref{thm:generalrank} bounds the number of halvings by $O(N^3)$ and the rounds per halving by $O(N^6)$. Once $g<1$, the split on $g$ enforces $g\le0$ in one further round. Both this cut and its split satisfy the length bound $p(L)$. Taking the maximum of these round counts over the defining inequalities, we obtain $R_r=P_I$ for some $r=O(N^9)$. If $P\ne\varnothing$ and $P_I=\varnothing$, one halving with $g\equiv1$ and $H=1$ forces $R_b=\varnothing$. If $P=\varnothing$, then $R_0=P_I$. Within each round, the finitely many eligible cuts can also be imposed successively: after any cut, both branches of each remaining split are subsets of their branches at the start of the round, so every remaining cut is still valid.
\end{proof}

\section*{AI use disclosure}
During the preparation of this manuscript, the authors used OpenAI's ChatGPT to find relevant literature outside of their knowledge, the development of induction arguments and checking the encoding bound, and the parts of the exposition including figures. The authors have verified and take full responsibility for the manuscript including the results, their proofs, and references.

\bibliographystyle{amsplain}
\bibliography{rank_references}
\end{document}